\documentclass[12pt]{article}
\usepackage{amssymb,amsmath}
\usepackage{hyperref}
\usepackage{graphicx}
\usepackage{color}

\begin{document}

\title{\LARGE\bf Identities for the arctangent function by enhanced midpoint integration and the high-accuracy computation of pi}

\author{
\normalsize\bf S. M. Abrarov\footnote{\scriptsize{Dept. Earth and Space Science and Engineering, York University, Toronto, Canada, M3J 1P3.}}\, and B. M. Quine$^{*}$\footnote{\scriptsize{Dept. Physics and Astronomy, York University, Toronto, Canada, M3J 1P3.}}}

\date{April 10, 2016}
\maketitle

\begin{abstract}
We describe a method of integration to obtain identities of the arctangent function and show how this method can be applied to the high-accuracy computation of the constant pi using the equation $\pi = 4 \arctan \left( 1 \right)$. Our approach combines the midpoint method with the Taylor expansion series to enhance accuracy in the subintervals. The accuracy of this method of integration is determined by number of subintervals $L$ and by order of the Taylor expansion $M$. This approach provides significant flexibility in computation since the required convergence in resulting equations can be optimized through appropriate choices for the integers $L$ and $M$. Sample computations are presented to illustrate that even with relatively small values of the integers $L$ and $M$ the constant $\pi$ can be computed with high accuracy.
\vspace{0.25cm}
\vspace{0.25cm}
\end{abstract}

\section{Introduction}

Consider a function $f\left( t \right)$ that is defined in the interval $t \in \left[ {a,b} \right]$. The integral of this function
$$
I = \int\limits_a^b {f\left( t \right)dt}
$$
can be computed numerically by using different efficient computational methods. In integration methods based on the Newton--Cotes formulas, such as the midpoint rule, the trapezoidal rule, the Simpson\text{'}s rule, the Boole\text{'}s rule and so on, the interval $\left[ {a, b} \right]$ is divided into $L$ sufficiently small subintervals  \cite{Kythe2005, Beu2015}
$$
\Delta {t_\ell } = {t_\ell } - {t_{\ell  - 1}}, \qquad \ell  = \left\{ {1,2,3, \,\, \ldots \,\, L} \right\},
$$
where
$$
t_0 = a, t_1 = a + \Delta t_1, t_2 = a + \Delta t_1 + \Delta t_2, \,\, \ldots \,\, t_L = b
$$
and
$$
\Delta {t_1} + \Delta {t_2} + \Delta {t_3} + \,\, \ldots \,\, + \Delta {t_L} = b - a.
$$

For convenience we assume that $a = 0$ and $b = 1$ since the described procedure can be easily generalized to any arbitrary values of $a$ and $b$. Consequently, we have
$$
I = \int\limits_0^1 {f\left( t \right)dt}.
$$
Implying that all subintervals $\Delta {t_\ell }$ are equal 
$$
\Delta {t_1} = \Delta {t_2} = \Delta {t_3} \ldots  = \Delta {t_L} = 1/L
$$
we can write
\begin{equation}
\label{eq_1}
I = \mathop {\lim }\limits_{\Delta {t_\ell } \to 0} \sum\limits_{\ell  = 1}^L {f\left( {{t_{\ell  - 1}} + \frac{{\Delta {t_\ell }}}{2}} \right)\Delta {t_\ell }} = \mathop {\lim }\limits_{L \to \infty } \frac{1}{L}\sum\limits_{\ell  = 1}^L {f\left( {\frac{{\ell  - 1/2}}{L}} \right)}.
\end{equation}
The truncated form of this equation is known as the midpoint rule \cite{Kythe2005, Beu2015}
\begin{equation}
\label{eq_2}
I = \int\limits_0^1 {f\left( t \right)dt} \approx \frac{1}{L}\sum\limits_{\ell  = 1}^L {f\left( {\frac{{\ell  - 1/2}}{L}} \right)},	\qquad L >  > 1.
\end{equation}

\section{Numerical integration method}

Although the midpoint rule is one of the simplest methods, it can provide, nevertheless, a relatively high accuracy in numerical integration. However, an improvement in accuracy with this method requires a large quantity of subintervals $\Delta {t_\ell }$ that may cause some side effects including long computational time, truncation and rounding errors.

This problem can be effectively resolved by improving accuracy for each subintervals $\Delta {t_\ell }$. In particular, in this work we show that one of the efficient ways to improve accuracy within these subintervals is to apply the Taylor expansion series
\begin{equation}
\label{eq_3}
f\left( t \right) = \sum\limits_{m = 0}^\infty  {\frac{{{{\left. {{f^{\left( m \right)}}\left( t \right)} \right|}_{t = {t_s}}}}}{{m!}}{{\left( {t - {t_s}} \right)}^m}},
\end{equation}
where the notation ${\left. {{f^{\left( m \right)}}\left( t \right)} \right|_{t = {t_s}}}$ signifies the ${m^{{\rm{th}}}}$ derivative of the function $f\left( t \right)$ at some point $t = {t_s}$. Thus, by using the Taylor expansion series \eqref{eq_3} at the midpoints over each subinterval 
$$
{t_{\ell  - 1}} + \frac{{\Delta {t_\ell }}}{2} = \frac{{\ell  - 1/2}}{L}
$$
we obtain exact values for all small areas $\Delta {I_\ell }$ under the curve of the function $f \left( t \right)$ corresponding to subintervals $\Delta {t_\ell } = 1/L  = \ell/L -\left(\ell - 1 \right)/L$ as follows
\begin{equation}
\label{eq_4}
\Delta {I_\ell } = \int\limits_{\left( {\ell  - 1} \right)/L}^{\ell /L} {\sum\limits_{m = 0}^\infty  {\frac{{{{\left. {{f^{\left( m \right)}}\left( t \right)} \right|}_{t = \frac{{\ell  - 1/2}}{L}}}}}{{m!}}{{\left( {t - \frac{{\ell  - 1/2}}{L}} \right)}^m}dt} }.
\end{equation}

The total area $I$ under the curve of the function $f\left( t \right)$ within interval $\left[ {0,1} \right]$ is a superposition of these small areas 
\begin{equation}
\label{eq_5}
I = \sum\limits_{\ell  = 1}^L {\Delta {I_\ell }}.
\end{equation}
Therefore, substituting equation \eqref{eq_4} into equation \eqref{eq_5} leads to
$$
I = \sum\limits_{\ell  = 1}^L {\int\limits_{\left( {\ell  - 1} \right)/L}^{\ell /L} {\sum\limits_{m = 0}^\infty  {\frac{{{{\left. {{f^{\left( m \right)}}\left( t \right)} \right|}_{t = \frac{{\ell  - 1/2}}{L}}}}}{{m!}}{{\left( {t - \frac{{\ell  - 1/2}}{L}} \right)}^m}dt} } }
$$
or
\begin{equation}
\label{eq_6}
I = \int\limits_0^1 {f\left( t \right)dt}  = \sum\limits_{\ell  = 1}^L {\sum\limits_{m = 0}^\infty  {\frac{{{{\left( { - 1} \right)}^m} + 1}}{{{{\left( {2L} \right)}^{m + 1}}\left( {m + 1} \right)!}}{{\left. {{f^{\left( m \right)}}\left( t \right)} \right|}_{t = \frac{{\ell  - 1/2}}{L}}}} }
\end{equation}
since
\footnotesize
$$
\int\limits_{\left( {\ell  - 1} \right)/L}^{\ell /L} {\frac{{{{\left. {{f^{\left( m \right)}}\left( t \right)} \right|}_{t = \frac{{\ell  - 1/2}}{L}}}}}{{m!}}{{\left( {t - \frac{{\ell  - 1/2}}{L}} \right)}^m}dt}  = \frac{{\left. {{f^{\left( m \right)}}\left( t \right)} \right|_{t = \frac{{\ell  - 1/2}}{L}}}}{{m!}}\int\limits_{\left( {\ell  - 1} \right)/L}^{\ell /L} {{{\left( {t - \frac{{\ell  - 1/2}}{L}} \right)}^m}dt},
$$
\normalsize
and
$$
\int\limits_{\left( {\ell  - 1} \right)/L}^{\ell /L} {{{\left( {t - \frac{{\ell  - 1/2}}{L}} \right)}^m}dt}  = \frac{{{{\left( { - 1} \right)}^m} + 1}}{{{{\left( {2L} \right)}^{m + 1}}\left( {m + 1} \right)}} = \frac{{\left({{\left( { - 1} \right)}^m} + 1 \right)m!}}{{{{\left( {2L} \right)}^{m + 1}}\left( {m + 1} \right)!}}.
$$

The equation \eqref{eq_6} is valid only if the Taylor expansion series \eqref{eq_4} over each midpoint remains non-truncated. However, if we attempt to truncate it up to an order $M$ as given by
\begin{equation}
\label{eq_7}
\Delta {I_\ell } \approx \int\limits_{\left( {\ell  - 1} \right)/L}^{\ell /L} {\sum\limits_{m = 0}^M {\frac{{{{\left. {{f^{\left( m \right)}}\left( t \right)} \right|}_{t = \frac{{\ell  - 1/2}}{L}}}}}{{m!}}{{\left( {t - \frac{{\ell  - 1/2}}{L}} \right)}^m}dt} }
\end{equation}
then instead of equation \eqref{eq_5} we have to use the limit
\begin{equation}
\label{eq_8}
I = \mathop {\lim }\limits_{L \to \infty } \sum\limits_{\ell  = 1}^L {\Delta {I_\ell }}
\end{equation}
since the accuracy of the equation \eqref{eq_7} tends to absolute only if $\Delta {t_\ell }$ tend to zero. Consequently, substituting the approximation \eqref{eq_7} into the limit \eqref{eq_8} leads to
$$
I = \mathop {\lim }\limits_{\ell  \to \infty } \sum\limits_{\ell  = 1}^L {\int\limits_{\left( {\ell  - 1} \right)/L}^{\ell /L} {\sum\limits_{m = 0}^M {\frac{{{{\left. {{f^{\left( m \right)}}\left( t \right)} \right|}_{t = \frac{{\ell  - 1/2}}{L}}}}}{{m!}}{{\left( {t - \frac{{\ell  - 1/2}}{L}} \right)}^m}dt} } }
$$
or
\begin{equation}
\label{eq_9}
I = \int\limits_0^1 {f\left( t \right)dt}  = \mathop {\lim }\limits_{L \to \infty } \sum\limits_{\ell  = 1}^L {\sum\limits_{m = 0}^M {\frac{{{{\left( { - 1} \right)}^m} + 1}}{{{{\left( {2L} \right)}^{m + 1}}\left( {m + 1} \right)!}}{{\left. {{f^{\left( m \right)}}\left( t \right)} \right|}_{t = \frac{{\ell  - 1/2}}{L}}}} }.
\end{equation}
It should be noted that the odd values $m$ do not contribute to integration since in equation \eqref{eq_9} the numerator ${\left( { - 1} \right)^m} + 1 = 0$ at $m = \left\{ {1,3,5, \ldots } \right\}.$

Is not difficult to see now that the conventional midpoint rule  \eqref{eq_2} is, in fact, the roughest method for approximation as it accounts only for the first term in the Taylor expansion series \eqref{eq_7} at the midpoints $\left( {\ell  - 1/2} \right)/L$. In particular, truncating the equation \eqref{eq_9} up to the order $M = 0$ we immediately obtain the limit \eqref{eq_1} corresponding to the conventional midpoint rule.

\section{Application to the arctangent function}

The arctangent function is one of the interesting functions that are widely used in many applications. Historically, many identities for the arctangent function were discovered centuries ago \cite{Beckmann1976, Berggren2004}. However, finding new equations for the arctangent function still remains topical in Applied Mathematics and Computational Physics, and several efficient identities for the arctangent function have been reported in the modern literature \cite{Calcut1999, Chen2010, Sofo2012}.

The arctangent function can be defined by the integral
\begin{equation}
\label{eq_10}
\arctan \left( x \right) = \int\limits_0^1 {\frac{x}{{1 + {x^2}{t^2}}}dt}.
\end{equation}
Consequently, substituting its integrand into equation \eqref{eq_9} yields
\small
\begin{equation}
\label{eq_11}
\arctan \left( x \right) = \mathop {\lim }\limits_{L \to \infty } \sum\limits_{\ell  = 1}^L {\sum\limits_{m = 0}^M {\frac{{{{\left( { - 1} \right)}^m} + 1}}{{{{\left( {2L} \right)}^{m + 1}}\left( {m + 1} \right)!}}\frac{{{\partial ^m}}}{{\partial {t^m}}}{{\left. {\left( {\frac{x}{{1 + {x^2}{t^2}}}} \right)} \right|}_{t = \frac{{\ell  - 1/2}}{L}}}} }.
\end{equation}
\normalsize

We consider only three examples $M = 0$, $M = 2$ and $M = 6$. At $M = 0$ the equation \eqref{eq_11} provides \cite{Abrarov2016}
\begin{equation}
\label{eq_12}
\arctan \left( x \right) = \mathop {\lim }\limits_{L \to \infty } \sum\limits_{\ell  = 1}^L {\frac{{4Lx}}{{4{L^2} + {{\left( {2\ell  - 1} \right)}^2}{x^2}}}}.
\end{equation}
This limit can also be obtained directly from the integral \eqref{eq_10} by using the equation \eqref{eq_1}. At $M = 2$ the equation \eqref{eq_11} results in
\small
\begin{equation}
\label{eq_13}
\arctan \left( x \right) = \mathop {\lim }\limits_{L \to \infty } \sum\limits_{\ell  = 1}^L {\left( {\frac{{4Lx}}{{4{L^2} + {{\left( {2\ell  - 1} \right)}^2}{x^2}}} - \frac{{4L{x^3}\left( {4{L^2} - 3{{\left( {2\ell  - 1} \right)}^2}{x^2}} \right)}}{{3{{\left( {4{L^2} + {{\left( {2\ell  - 1} \right)}^2}{x^2}} \right)}^3}}}} \right)}.
\end{equation}
\normalsize
while at $M = 6$ it follows that
\footnotesize
\begin{equation}
\label{eq_14}
\begin{aligned}
\arctan & \left( x \right) = \mathop {\lim }\limits_{L \to \infty } \sum\limits_{\ell  = 1}^L {\left( {\frac{{4Lx}}{{4{L^2} + {{\left( {2\ell  - 1} \right)}^2}{x^2}}} - \frac{{4L{x^3}\left( {4{L^2} - 3{{\left( {2\ell  - 1} \right)}^2}{x^2}} \right)}}{{3{{\left( {4{L^2} + {{\left( {2\ell  - 1} \right)}^2}{x^2}} \right)}^3}}}} \right.} \\
& + \frac{{4L{x^5}\left( {16{L^4} - 40{{\left( {2\ell  - 1} \right)}^2}{L^2}{x^2} + 5{{\left( {2\ell  - 1} \right)}^4}{x^4}} \right)}}{{5{{\left( {4{L^2} + {{\left( {2\ell  - 1} \right)}^2}{x^2}} \right)}^5}}} \\
& \left. { - \frac{{4L{x^7}\left( {64{L^6} - 336{{\left( {2\ell  - 1} \right)}^2}{L^4}{x^2} + 140{{\left( {2\ell  - 1} \right)}^4}{L^2}{x^4} - 7{{\left( {2\ell  - 1} \right)}^6}{x^6}} \right)}}{{7{{\left( {4{L^2} + {{\left( {2\ell  - 1} \right)}^2}{x^2}} \right)}^7}}}} \right).
\end{aligned}
\end{equation}
\normalsize

As the equation \eqref{eq_11} consists of double summation with respect to indices $\ell $ and $m$, its algorithmic implementation may be advantageous in optimization procedure. In particular, the equation \eqref{eq_11} can provide some flexibility to find a best combination of the upper values  $L$ and $M$ in truncation to perform more efficient numerical integration. Furthermore, the equation \eqref{eq_11} leads to infinite quantity of identities for the arctangent function because the upper integer $M$ may be taken arbitrarily large.

\section{Examples of pi identities and computation}

There are many identities for the constant pi expressed in terms of the arctangent function that can be written in form \cite{Abeles1993, Lehmer1938, Borwein1989, Borwein2015}
$$
\pi  = \sum\limits_{n = 1}^N {{\alpha _n}\arctan \left( {\frac{1}{{{\beta _n}}}} \right)}  = \sum\limits_{n = 1}^N {{\alpha _n}{\rm{arccot}}\left( {{\beta _n}} \right)},
$$ 
where ${\alpha _n}$ and ${\beta _n}$ are some integers. However, we take the simplest identity to compute pi
\begin{equation}
\label{eq_15}
\pi  = 4\arctan \left( 1 \right).
\end{equation}
Thus, substituting $x = 1$ into equations \eqref{eq_12}, \eqref{eq_13}, \eqref{eq_14} and using identity \eqref{eq_15} we obtain the following limits for pi
\begin{equation}
\label{eq_16}
\pi  = 4 \times \mathop {\lim }\limits_{L \to \infty } \sum\limits_{\ell  = 1}^L {\frac{{4L}}{{4{L^2} + {{\left( {2\ell  - 1} \right)}^2}}}},
\end{equation}
\begin{equation}
\label{eq_17}
\pi  = 4 \times \mathop {\lim }\limits_{L \to \infty } \sum\limits_{\ell  = 1}^L {\left( {\frac{{4L}}{{4{L^2} + {{\left( {2\ell  - 1} \right)}^2}}} - \frac{{4L\left( {4{L^2} - 3{{\left( {2\ell  - 1} \right)}^2}} \right)}}{{3{{\left( {4{L^2} + {{\left( {2\ell  - 1} \right)}^2}} \right)}^3}}}} \right)}
\end{equation}
and
\small
\begin{equation}
\label{eq_18}
\begin{aligned}
\pi  =& \, 4 \times \mathop {\lim }\limits_{L \to \infty } \sum\limits_{\ell  = 1}^L {\left( {\frac{{4L}}{{4{L^2} + {{\left( {2\ell  - 1} \right)}^2}}} - \frac{{4L\left( {4{L^2} - 3{{\left( {2\ell  - 1} \right)}^2}} \right)}}{{3{{\left( {4{L^2} + {{\left( {2\ell  - 1} \right)}^2}} \right)}^3}}}} \right.} \\
 &+ \frac{{4L\left( {16{L^4} - 40{{\left( {2\ell  - 1} \right)}^2}{L^2} + 5{{\left( {2\ell  - 1} \right)}^4}} \right)}}{{5{{\left( {4{L^2} + {{\left( {2\ell  - 1} \right)}^2}} \right)}^5}}}\\
&\left. { - \frac{{4L\left( {64{L^6} - 336{{\left( {2\ell  - 1} \right)}^2}{L^4} + 140{{\left( {2\ell  - 1} \right)}^4}{L^2} - 7{{\left( {2\ell  - 1} \right)}^6}} \right)}}{{7{{\left( {4{L^2} + {{\left( {2\ell  - 1} \right)}^2}} \right)}^7}}}} \right).
\end{aligned}
\end{equation}
\normalsize

In order to compare accuracy and convergence rates in these equations, we performed sample computations by using Wolfram Mathematica 9 in enhanced precision mode. The sample computations are performed at same value of the integer $L = 1000$ in all three equations above.

The actual value of the constant $\pi $ is given by
$$
3.1415926535897932384626433832795028841971693993751 \ldots
$$
while the results of computation provided by equations \eqref{eq_16}, \eqref{eq_17} and \eqref{eq_18} are
$$
\underbrace {3.141592}_{7\,\,{\rm{digits}}}7369231265717940545935969641467776336373917 \ldots \,\,,
$$
$$
\underbrace {3.14159265358979323846}_{21\,\,{\rm{digits}}}37594547080737075957760027852 \ldots 
$$
and
$$
\underbrace {3.1415926535897932384626433832795028}_{35\,\,{\rm{digits}}}649618474297397\ldots 
$$
respectively. Comparing these results with actual value for the constant $\pi $ we can see a rapid improvement in convergence with increasing $M$. Specifically, at $M = 0$, $M = 2$ and $M = 6$ the quantities of coinciding digits are $7$, $21$ and $35$, respectively. Thus, with only $L = 1000$ and $M = 6$ the equation \eqref{eq_18} alone provides high-accuracy computation with $35$ coinciding digits. Observing such a tendency we may expect that further increase of the integer $M$ in the equation \eqref{eq_11} would greatly improve the convergence in computation of $\pi$.

Due to simple basic rules that are used for derivatives, taking the $m^{\rm{th}}$ derivative of the function ${\frac{x}{1+x^2t^2}}$ at the midpoints ${t = \frac{{\ell  - 1/2}}{L}}$ is not problematic even if $m > > 1$. Consequently, assuming $x = 1$ it is not difficult to find algorithmically the expansion coefficients
$$
g_{\ell, m} = \frac{{{d ^m}}}{{d {t^m}}}{{\left. {\left( {\frac{1}{{1 + {t^2}}}} \right)} \right|}_{t = \frac{{\ell  - 1/2}}{L}}}
$$
for the equations
\[
\pi = 4 \times \mathop {\lim }\limits_{L \to \infty } \sum\limits_{\ell  = 1}^L {\sum\limits_{m = 0}^M} \frac{{{{\left( { - 1} \right)}^m} + 1}}{{{{\left( {2L} \right)}^{m + 1}}\left( {m + 1} \right)!}} g_{\ell,m}
\]
and
\[
\pi = 4 \times \sum\limits_{\ell  = 1}^L {\sum\limits_{m = 0}^\infty} \frac{{{{\left( { - 1} \right)}^m} + 1}}{{{{\left( {2L} \right)}^{m + 1}}\left( {m + 1} \right)!}} g_{\ell,m}.
\]
Thus, the computational test we performed shows that at $M = L = 46$ the number of coinciding digits reaches $105$. Since the odd values of $m$ do not contribute, the total number of the terms $\frac{{{{\left( { - 1} \right)}^m} + 1}}{{{{\left( {2L} \right)}^{m + 1}}\left( {m + 1} \right)!}}g_{\ell, m}$ involved in summation is $L \times \left( M/2 + 1 \right) = 46 \times \left( 46/2 + 1 \right) = 1104$.

\section{Conclusion}

We consider a numerical integration method where accuracy in subintervals can be improved by using the Taylor expansion series over the midpoints and show some identities for the arctangent function as examples. Since the accuracy of this method of integration is determined by number of subintervals $L$ and by order of the Taylor expansion $M$, a required convergence in obtained equations of the arctangent function can be optimized by choosing a combination of the integers $L$ and $M$. Sample computations reveal that this approach provides a high-accuracy in computation of the constant pi even at relatively small values of the integers $L$ and $M$.

\section*{Acknowledgments}

This work is supported by National Research Council Canada, Thoth Technology Inc. and York University.



\end{document}